\documentclass{elsart3-1}
\usepackage{amsmath,amssymb,amsfonts}
\usepackage{mathabx,stmaryrd}
\usepackage{bbold}
\usepackage{paralist}
\usepackage[normalem]{ulem}
\usepackage{nicefrac}
\usepackage{xspace}
\usepackage{xcolor}
\usepackage{tikz}
\usepackage{caption}
\usepackage{subcaption}
\usepackage{pgfplots,pgfplotstable}
\pgfqkeys{/pgfplots}{
  cycle list name = black white
}
\pgfplotsset{compat=newest}

\usepackage[english,francais]{babel}

\renewcommand{\vec}[1]{\boldsymbol{#1}}

\renewcommand{\div}{{\rm div}}

\newcommand{\R}{\mathbb{R}}
\newcommand{\E}{{\mathbb E}}

\newcommand{\bigo}{{\mathcal O}}

\DeclareMathOperator{\ddiv}{div}
\DeclareMathOperator{\trace}{trace}

\newcommand{\cA}{\mathcal A} 
\newcommand{\cB}{\mathcal B} 
\newcommand{\cL}{\mathcal L} 
\newcommand{\cH}{\mathcal H} 

\newtheorem{theorem}{Theorem}[section]

\newtheorem{e-proposition}[theorem]{Proposition}

\newtheorem{e-definition}[theorem]{Definition\rm}
\newtheorem{remark}{\it Remark\/}

\renewcommand{\theequation}{\arabic{equation}}
\newenvironment{psmallmatrix}{\left(\begin{smallmatrix}}{\end{smallmatrix}\right)}

\def\og{\leavevmode\raise.3ex\hbox{$\scriptscriptstyle\langle\!\langle$~}}
\def\fg{\leavevmode\raise.3ex\hbox{~$\!\scriptscriptstyle\,\rangle\!\rangle$}}

\journal{the Acad\'emie des Sciences}
\begin{document}
\centerline{}
\begin{frontmatter}


\selectlanguage{english}
\vskip-40mm
\title{Accelerated convergence to equilibrium and reduced asymptotic variance for Langevin dynamics using Stratonovich perturbations}


\selectlanguage{english}
\author[1]{Assyr Abdulle\thanksref{aa}},
\ead{Assyr.Abdulle@epfl.ch}
\author[2]{Grigorios A. Pavliotis\thanksref{gp}},
\ead{g.pavliotis@imperial.ac.uk}
\author[3]{Gilles Vilmart\thanksref{gv}}
\ead{Gilles.Vilmart@unige.ch}

\thanks[aa]{Partially supported by the Swiss National Science Foundation, grant No. 200020\_172710.}

\thanks[gp]{Partially supported by the EPSRC, grants No. EP/P031587/1 	 EP/L024926/1 	EP/L020564/1.}

\thanks[gv]{Partially supported by the Swiss National Science Foundation, grants No. 200021\_162404, No. 200020\_178752 and No. 200020\_144313/1.\\
The authors thank Andrew Duncan and Tony Leli\`evre for useful discussions. 
They acknowledge the kind hospitality of the Erwin Schr\"odinger International Institute for Mathematics and Physics (ESI), where part of this research was developed under the frame of the Thematic Programme {\it Numerical Analysis of Complex PDE Models in the Sciences}. Part of the work was done while GP was visiting the Chair of Computational Mathematics and Numerical Analysis at EPFL. The hospitality and financial support is greatly acknowledged.}

\address[1]{\'Ecole Polytechnique F\'ed\'erale de Lausanne (EPFL), SB-MATH-ANMC,
          Station 8,
          1015 Lausanne, Switzerland}
          \address[2]{Department of Mathematics, Imperial College London, London SW7 2AZ, UK}
          \address[3]{Universit\'e de Gen\`eve, Section de math\'ematiques, 2-4 rue du Li\`evre, CP 64, CH-1211 Gen\`eve 4, Switzerland}


\vskip-0.9cm

\begin{abstract}
\selectlanguage{english}
In this paper we propose a new approach for sampling from probability measures in, possibly, high dimensional spaces. By perturbing the standard overdamped Langevin dynamics by a suitable  Stratonovich perturbation that preserves the invariant measure of the original system, we show that accelerated convergence to equilibrium and reduced asymptotic variance can be achieved, leading, thus, to a computationally advantageous sampling algorithm. The new perturbed Langevin dynamics is reversible with respect to the target probability measure and, consequently, does not suffer from the drawbacks of the nonreversible Langevin samplers that were introduced in~[C.-R. Hwang, S.-Y. Hwang-Ma, and S.-J. Sheu, Ann. Appl. Probab. 1993] and studied in, e.g. [T. Lelievre, F. Nier, and G. A. Pavliotis J. Stat. Phys., 2013] and [A. B. Duncan, T. Leli\`evre, and G. A. Pavliotis J. Stat. Phys., 2016], while retaining all of their advantages in terms of accelerated convergence and reduced asymptotic variance. In particular, the reversibility of the dynamics ensures that there is no oscillatory transient behaviour. The improved performance of the proposed methodology, in comparison to the standard overdamped Langevin dynamics and its nonreversible perturbation, is illustrated on an example of sampling from a two-dimensional warped Gaussian target distribution.

\vskip 0.5\baselineskip

\selectlanguage{french}
\noindent{\bf R\'esum\'e} \vskip 0.5\baselineskip \noindent
{\bf Convergence acc\'el\'er\'ee vers l'\'equilibre et variance asymptotique r\'eduite pour la dynamique de Langevin \`a l'aide de perturbations de Stratonovich.}
Dans cet article, nous proposons une nouvelle approche pour l'\'echantillonnage de mesures invariantes dans des espaces de grandes dimensions  \`a l'aide d'une dynamique de Langevin perturb\'ee.
En modifiant la dynamique standard de l'\'equation de Langevin suramortie en introduisant une perturbation de Stratonovich convenable pr\'eservant la mesure invariante du syst\`eme initial, nous montrons qu'il est possible d'obtenir une convergence acc\'el\'er\'ee vers l'\'equilibre et une variance asymptotique r\'eduite, conduisant ainsi \`a un algorithme d'\'echantillonnage avantageux du point de vue du calcul. La nouvelle dynamique de Langevin perturb\'ee est r\'eversible par rapport \`a la mesure de probabilit\'e cherch\'ee et ne souffre donc pas des inconv\'enients des \'echantillonneurs de Langevin non r\'eversibles introduits dans [C.-R. Hwang, S.-Y. Hwang-Ma, and S.-J. Sheu, Ann. Appl. Probab. 1993]   et \'etudi\'es par exemple dans [T. Lelievre, F. Nier, and G. A. Pavliotis J. Stat. Phys., 2013] et [A. B. Duncan, T. Leli\`evre, and G. A. Pavliotis J. Stat. Phys., 2016], tout en conservant tous leurs avantages en termes de convergence acc\'el\'er\'ee et de r\'eduction de la variance asymptotique. En particulier, la r\'eversibilit\'e de la dynamique garantit l'absence de comportement transitoire oscillant. Les performances am\'elior\'ees de la m\'ethodologie propos\'ee, par rapport \`a la dynamique de Langevin suramortie standard et \`a sa perturbation irr\'eversible, sont illustr\'ees par un exemple d'\'echantillonnage \`a partir d'une distribution  gaussiennes d\'eform\'ees \`a deux dimensions.

\end{abstract}
\end{frontmatter}

\selectlanguage{english}
\section{Introduction} \label{se:intro}
Sampling from probability measures in high dimensional spaces is an important problem that arises in several applications, including computational statistical physics~\cite{Lelievre_al_2009}, Bayesian inference~\cite{Stuart2010}, and machine learning~\cite{Andrieu2003}. Typically one is interested in calculating integrals of the form
\begin{equation}\label{e:expect}
\pi(\phi):= \E_{\pi} \phi : = \int_{\R^d}\phi(x) \, \pi(dx),
\end{equation}
where $\pi(dx)  =\pi(x) \, dx$\footnote{We assume that the target probability measure has a density with respect to Lebesgue measure. To simplify the notation, we will denote both the measure and the density by $\pi$.} is a probability measure in $\R^d$, known up to the normalization constant and $\phi \in L^2(\pi)$ is an observable. Here  $L^2(\pi)$ denotes  the weighted $L^2$ space for the scalar  product $(\phi,\psi)_{\pi}=\int_{\R^d} \phi(x)\psi(x)\pi(x) dx$  and the corresponding norm is denoted by $\|\phi\|_{\pi}$.
A standard methodology for calculating, or, rather, estimating the integral in~\eqref{e:expect} is to construct a stochastic process $\{ X(t) \}_{t >0}$ in $\R^d$, e.g. an It\^{o} diffusion process 
\begin{equation} \label{eq:sde0}
dX(t)= f(X(t)) \, dt + \sigma(X(t)) \, dW_t
\end{equation}
that is ergodic with respect to the measure $\pi$. Here $W_t$ is a standard $m$--dimensional Brownian motion and  $f:\R^d\rightarrow \R^d$ and $\sigma:\R^d\rightarrow \R^{d\times m}$ are assumed smooth and Lipschitz continuous. In particular, $\pi$ is the unique normalized solution of the stationary Fokker-Plank equation $\cL^*\pi=0,$ where 
$\cL^*$ is the $L^2(dx)$ adjoint of the generator  $\cL\phi:=f\cdot\nabla\phi + \frac{1}{2}\sigma\sigma^T: \nabla^2\phi$ of the SDE \eqref{eq:sde0}.\footnote{For two matrices $A$ and $B$ we use the notation $A:B=\trace(A^TB)$.}  In what follows we  denote
by $\cH^*$ the $L^2(dx)$ adjoint of an operator $\cH$ and by  $\cH^\sharp$ its $L^2(\pi)$ adjoint.

Under appropriate assumptions on the drift and diffusion coefficients, we can prove a strong law of large numbers and a central limit theorem as  $T \rightarrow \infty,$
\begin{equation}\label{e:time_conv}
\pi_T(\phi):= \frac{1}{T} \int_0^T \phi(X(t)) \, dt  \rightarrow \pi(\phi) \quad \mbox{a.e.}, \; \; X_0 =  x,
\end{equation}
and we have the following convergence in law 
\begin{equation}\label{e:clt}
\sqrt{T} (\pi_T(\phi) - \pi(\phi)) \rightarrow \mathcal{N} (0, \sigma_{\phi}^2),
\end{equation}
where $\sigma^2_{\phi}$ denotes the asymptotic variance of the observable $\phi$,  
 given by the Kipnis-Varadhan formula
\begin{equation}\label{e:kipnis}
\sigma_{\phi}^2 = \langle (\phi - \pi (\phi), (-\cL)^{-1}(\phi - \pi (\phi)) \rangle_{\pi}.
\end{equation}
Under the assumption that the generator has a spectral gap in $L^2(\pi)$ (see for instance \cite{Mattingly10con}) we have the following exponential convergence
\begin{equation}\label{e:exp}
\big|\E(\phi(X(t)))-\pi(\phi)\big|\leq Ce^{-\lambda t},
\end{equation}
where $\lambda>0$ is the spectral gap of the generator $\cL$.

In this paper, we focus on the overdamped Langevin dynamics for sampling \eqref{e:expect},
\begin{equation}\label{e:langevin}
d X(t) = f(X(t)) \, dt + \sqrt{2} \, dW_t,
\end{equation}
where $f(x):=-\nabla V(x)$ and $W_t$ is a standard $d$-dimensional Brownian motion.
The invariant measure of~\eqref{1} is given by $\pi(dx)=Z^{-1}e^{-V(x)} \, dx,$ where $Z = \int_{\R^d} e^{-V(x)} \, dx$ is the normalization constant and $V:\R^d\rightarrow \R$ is a smooth confining potential.
A question that has attracted considerable attention in recent years is the construction of modified Langevin dynamics that have better sampling properties in comparison to the standard overdamped Langevin dynamics \eqref{e:langevin}.  
Several modifications of the Langevin dynamics~\eqref{e:langevin} that can be used in order to sample from $\pi$ are presented in~\cite[Sec 2.2]{DuncanNuskenPavliotis2017}. A well known technique that was first  introduced in~\cite{Hwang_al1993,Hwang_al2005} and analyzed in a series of recent papers, e.g.~\cite{Bellet_Spiliopoulos_2015,Bellet_Spiliopoulos_2016,LelievreNierPavliotis2013,DuncanLelievrePavliotis2016}
for improving the performance of the Langevin sampler~\eqref{e:langevin}, is to introduce in \eqref{eq:sde0} a divergence-free (with respect to the target distribution) drift perturbation $g:\R^d\rightarrow \R^d,$
\begin{equation} \label{eq:sdedet}
dX(t) = (f(X(t)) + g(X(t))) dt + \sqrt{2} \,dW_t, 
\end{equation}
such that
\begin{equation}\label{eq:divg}
\ddiv(g\pi)=0.
\end{equation}
We will refer to~\eqref{eq:divg} as the divergence-free condition. This condition ensures that the SDE \eqref{eq:sdedet} has the same invariant  measure $\pi$ as \eqref{eq:sde0}. 
We remark that there are infinitely many vector fields $g$ that satisfy~\eqref{eq:divg}. A complete characterization of all vector fields that satisfy this condition can be found in~\cite[Prop. 2.2]{Hwang_al2005}. 

It is by now a standard, and not difficult to prove, result that nonreversible dynamics exhibits better properties as a sampling scheme, in the sense that the nonreversible perturbation accelerates convergence to equilibrium and reduces the asymptotic variance. The generator of the nonreversible dynamics \eqref{eq:sdedet} is given by 
\begin{equation} \label{eq:defLD}
\cL_D  = \cL  + \cA,
\end{equation}
where $\cL$ is the generator of \eqref{eq:sde0} $\cA$ is defined by $\cA\phi=\nabla\phi\cdot g$ (in the calculations below we will use the notation $\cA\phi=\phi'g$). 
The drawback of the nonreversible Langevin sampler~\eqref{eq:sdedet} is that, since the generator of the dynamics is a nonselfadjoint operator, a transient, oscillatory phase is introduced.
This transient behaviour can be addressed, in principle, by the use of an appropriate splitting numerical scheme~\cite{DuncanPavliotisZygalakis2017}. 

In this paper, we introduce and analyze an alternative way for perturbing the overdamped reversible Langevin dynamics that is reversible
and enjoys all the advantages of the nonreversible sampler \eqref{eq:sdedet}, while not suffering from the drawback of its oscillatory transient dynamics. 
The new dynamics is given by the Stratonovitch perturbation \begin{equation} \label{eq:sdestrato}
dX(t) = f(X(t))\,dt+ g(X(t)) \circ \sqrt{2}\,d\beta_t + \sqrt{2} \, dW_t,
\end{equation}
where $g$ satisfies the divergence-free condition and we assume that $\beta_t$ is a one-dimensional standard Wiener process that is independent of $W_t$.\footnote{One can also consider Stratonovich perturbations driven by multidimensional Brownian motions with diffusion functions $g^j,$ $j=1,2,\dots$ satisfying
$\ddiv(g^j\pi)=0.$ A detailed analysis of such perturbed Langevin dynamics will be presented elsewhere.}
For the Stratonovich-perturbed Langevin dynamics~\eqref{eq:sdestrato} we have the following result.
\begin{theorem}[Reversibility of the perturbed dynamics]
\label{th:reversibility}
Considered the perturbed dynamics \eqref{eq:sdestrato}, were $g$ satisfies the divergence-free condition~\eqref{eq:divg}. Then  the generator of  \eqref{eq:sdestrato} can be written in the form
\begin{equation}\label{eq:defLS}
\cL_S  = \cL  + \cA^2,
\end{equation}
and $\cL_S$ is symmetric in  $L^2(\pi)$, i.e. $\cL_S=\cL_S^\sharp$.
\end{theorem}
As a consequence of Theorem~\ref{th:reversibility}, the eigenvalues of $\cL_S$ are real, hence there is no transient behaviour of the dynamics. 
\begin{remark}
The proposed modified Langevin sampler~\eqref{eq:sdestrato} can be written in the form of general reversible diffusion process\footnote{Our results can be extended to cover the case of the preconditioned/Riemannian manifold Markov Chain Monte Carlo Langevin dynamics. The details will be presented elsewhere.} 
(see~\cite[Ch. 4]{Pavl2014} for the characterization of diffusion processes that are reversible with respect to a given measure):
$$
dX(t) = - (M\nabla V)(X(t)) \, dt + (\div M)(X(t)) \, dt + \sqrt{2} D(X(t)) \, d\widehat{W}_t, 
$$ 
where $M = I + g g^T=D D^T \in R^{d \times d}, \, D = (I,g) \in \R^{d \times (d+1)}$ and $\widehat{W}_t = (W_t, \beta_t)$ is a standard $d+1$ dimensional Brownian motion.
\end{remark}
\begin{theorem}[Invariant measure preservation]
\label{th:pres_inv_measure}
Under the assumptions of Theorem \ref{th:reversibility}, the perturbed dynamics~\eqref{eq:sdestrato} is ergodic with respect to the measure $\pi(dx) = Z^{-1} e^{-V} \, dx$.
\end{theorem}
\medskip
\begin{remark}
We note that Theorems \ref{th:reversibility} and \ref{th:pres_inv_measure} remain true for general ergodic SDEs \eqref{eq:sde0} with a Stratonovich perturbation,
\begin{equation} 
dX(t) =  f(X(t))\,dt+ g(X(t)) \circ \sqrt{2}\,d\beta_t + \sigma(X(t))dW_t, 
\end{equation}
where $g$ is a divergence-free vector field with respect to $\pi$ and $f:\R^d\rightarrow \R^d$ does not have a gradient structure. This includes, in particular, degenerate diffusions (e.g. when the diffusion matrix $\sigma \sigma^T$ is only positive semidefinite), for example the underdamped Langevin dynamics. Indeed the gradient structure is not used in the proofs. Note however that in the case where the functional form of $\pi$ is not explicitly known, it can be difficult to compute such a vector field.
\end{remark}
The next theorem shows that, in comparison to the original overdamped Langevin dynamics~\eqref{e:langevin}, the Stratonovich perturbation yields a larger spectral gap and a reduced asymptotic variance. Similarly to the nonreversible deterministic perturbation \eqref{eq:sdedet}, this hence leads to an improved reversible sampler for the invariant measure \eqref{e:expect}, both in terms of speeding up the convergence to equilibrium \eqref{e:exp} as well as in terms of reducing the asymptotic variance \eqref{e:kipnis}. When combined, these results provide us with improved performance when measured in the mean-squared error; see \cite[Sec.2.3]{DuncanNuskenPavliotis2017}.

We recall that, under the assumption that the potential $V$ grows sufficiently fast at infinity, the generator of both the standard Langevin and of the Stratonovich-perturbed dynamics have purely discrete spectrum. 
\begin{theorem}[Accelerated convergence and reduced asymptotic variance]
\label{th:reveribility} Let the assumption of Theorem \ref{th:reversibility} hold and 
let $\lambda_{L}$ and $\lambda_S$ denote the spectral gaps of the overdamped Langevin~\eqref{e:langevin} and of the Stratonovich-perturbed dynamics~\eqref{eq:sdestrato}, respectively. Then
\begin{equation}\label{e:eigs}
\lambda_L \leq \lambda_S.
\end{equation}
Let, furthermore $\phi \in L^2(\pi)$ and denote the corresponding asymptotic variances by $\sigma^2_L(\phi)$ and $\sigma^2_S(\phi)$. Then
\begin{equation}\label{e:var}
\sigma^2_L(\phi) \geq \sigma^2_S(\phi). 
\end{equation}
\end{theorem}
\begin{remark} 
When the target distribution is Gaussian, in particular for the two dimensional quadratic potential $V(x)=\frac12(x_1^2+\lambda x_2^2)$ with $\lambda \ll 1$, the standard Langevin dynamics \eqref{e:langevin} converges to equilibrium at the very slow rate $\lambda_L=\lambda$, and it was shown in~\cite{LelievreNierPavliotis2013} that a perturbation of the form
\begin{equation} \label{eq:gdelta}
g(x)=\delta^\theta J\nabla V(x), \qquad J=\begin{psmallmatrix} 0& 1 \\ -1 & 0 \end{psmallmatrix},
\end{equation}
with size $\delta \sim \lambda^{-1/2}$ and $\theta=1$ yields in the case of a nonreversible perturbation \eqref{eq:sdedet} an optimally improved convergence rate $\lambda_D=\bigo(1)$. For isotropic Gaussians, the optimally reduced asymptotic variance using a nonreversible perturbation can also be calculated~\cite[Sec. 4]{DuncanLelievrePavliotis2016}. 
Similarly, an improved convergence rate of $\lambda_S=\bigo(1)$ can also be obtained for the reversible perturbation \eqref{eq:sdestrato} for the same scaling $\delta \sim \lambda^{-1/2}$ and $\theta=1/2$. 
Observe that the factor $\delta^\theta$ in \eqref{eq:gdelta} yields a perturbation of size $\mathcal{O}(\delta)$ of the Langevin generator $\cL$ in both perturbed generators $\cL_D$ in \eqref{eq:defLD} and $\cL_S$ in \eqref{eq:defLS}.
It is important to note that the optimal nonreversible perturbation depends on the optimality criterion used, i.e. on whether our aim is to maximize the rate of convergence to equilibrium or to minimize the asymptotic variance (uniformly over the space of square integrable observables). Contrary to this, the optimal reversible perturbation is the same with respect to these two optimality criteria. This observation will be explored further in a future work together with a complete analysis of optimal Stratonovitch perturbations for Gaussian target distributions.
\end{remark}

\section{Proof of the main results}

We start by recalling from~\cite{LelievreNierPavliotis2013} that the differential operator $\cA$ is skew-symmetric in  $L^2(\pi)$, i.e. $\cA^\sharp=-\cA$. This result follows from an integration by parts and \eqref{eq:divg}.
To prove our main results we also use that the original SDE \eqref{eq:sde0} has the the generator 
\begin{equation} \label{eq:defL}
\cL \phi = \phi' f +  \Delta\phi.
\end{equation}
{\it Proof of Theorem \ref{th:reversibility}}~
We convert the Stratonovitch SDE~\eqref{eq:sdestrato} into an It\^o one: 
\begin{equation} \label{eq:stratotointo}
dX = f(X) dt + g'(X)g(X) dt + g(X)  \sqrt{2}\,d\beta_t + \sqrt{2}\, dW_t.
\end{equation}
Using the calculation
$$
\cA^2\phi = (\phi'g)'g = \phi'g'g +  \phi''(g,g).
$$
we deduce the result \eqref{eq:defLS} by applying formula \eqref{eq:defL} to the SDE \eqref{eq:stratotointo}. An immediate consequence of 
$\cA^\sharp=-\cA$
is then that $(\cA^2)^\sharp=\cA^2$, i.e. $\cA^2$ is $L^2(\pi)$ symmetric. As $\cL$ itself is $L^2(\pi)$ symmetric, we have that $\cL_S$ is also $L^2(\pi)$ symmetric.
$\qed$
\smallskip\\
{\it Proof of Theorem \ref{th:pres_inv_measure}}
The $L^2$-adjoint satisfies
\begin{equation} 
\cL_S^* \pi = \cL^* \pi + \cA^{*} (\cA^{*} \pi)=0,
\end{equation}
where we have used the fact that $\cL^* \pi=\cA^{*} \pi =0$. Hence $\pi$ is the unique invariant measure of the perturbed dynamics \eqref{eq:sdestrato}.
$\qed$

\noindent
{\it Proof of Theorem \ref{th:reveribility}} 
We write the generator of the Stratonovich-perturbed dynamics as $\cL_S = - \cB^\sharp \cB - \cA^\sharp \cA$ with $\cB = \nabla, \;\cA = g\cdot \nabla, \; \cA^\sharp = - \cA $. The quadratic form associated to $\cL_S$ is $\langle - \cL_S \phi , \phi \rangle_{\pi} = \|\cB \phi \|_{\pi}^2 +  \|\cA \phi \|_{\pi}^2$ for all $\phi \in H^1(\pi)$ the weighted Sobolev space that is defined in the standard manner.  The quadratic form associated to the generator of the reversible Langevin dynamics $\cL = - \cB^\sharp \cB$ is $\langle - \cL \phi , \phi \rangle_{\pi} = \|\cB \phi \|_{\pi}^2$. Since both $\cL_S$ and $\cL$ are symmetric operators in $L^2(\pi)$ with compact resolvents, the spectral gap of the reversible Langevin dynamics is given by the Rayleigh quotient formula, 
\begin{eqnarray*}
\lambda_{S} = \min_{\phi \in H^1(\pi), \int \phi \pi = 0} \frac{\langle -\cL_S \phi , \phi \rangle_{\pi} }{\|\phi \|^2_{\pi}}  = \min_{\phi \in H^1(\pi), \int \phi \pi = 0} \frac{\|\cB \phi \|_{\pi}^2 +  \|\cA \phi \|_{\pi}^2 }{\|\phi \|^2_{\pi}} 
\geq  \min_{\phi \in H^1(\pi), \int \phi \pi = 0} \frac{\|\cB \phi \|_{\pi}^2  }{\|\phi \|^2_{\pi}} = \lambda_L.
\end{eqnarray*}

To prove the bound on the asymptotic variance, we first write the formula for $\sigma^2_S(\phi)$ in the form $\sigma^2_S(\phi) = \langle \psi_S,\phi  \rangle_{\pi}$ where $\psi_S$ is the solution of the Poisson equation $-\cL_S \psi_S = \phi,$ and where without loss of generality we have assumed that $\int_{\R^d} \phi \, \pi = 0$. 
We also consider $\psi_L$, the solution of the Poisson equation $-\cL \psi_L= \phi$ and using
$\cL=\cL_S+\cA^\sharp\cA$, we
obtain
\begin{eqnarray*}
\sigma^2_L(\phi) &=& \langle \phi ,\psi_L \rangle_{\pi} = \langle (- \cL_S)\psi_S,\psi_L \rangle_{\pi} = \langle\psi_S, (- \cL )\psi_L \rangle_{\pi} - \langle \cA^2\psi_S,\psi_L \rangle_{\pi} = \langle\psi_S, \phi \rangle_{\pi} + \langle \cA^\sharp\cA\psi_S,\psi_L \rangle_{\pi}\\
&=& \sigma^2_S(\phi)  + \langle \cA\psi_S, \cA\psi_L \rangle_{\pi}.
\end{eqnarray*}
To prove our claim, it is sufficient to show that $\langle \cA\psi_S, \cA\psi_L \rangle_{\pi} \geq 0$. We calculate,
\begin{eqnarray*}
\langle \cA\psi_S, \cA\psi_L \rangle_{\pi} & = & \langle \cA\psi_S, \cA (- \cL)^{-1} \phi \rangle_{\pi} = \langle \cA\psi_S, \cA (- \cL)^{-1} ((- \cL) +  (- \cA^2))\psi_S \rangle_{\pi}
\\ & = & \langle \cA^\sharp \cA\psi_S, (I +  (- \cL)^{-1} (- \cA^2))\psi_S \rangle_{\pi} = \|\cA\psi_S \|_{\pi}^2 +  \langle (- \cA^2)\psi_S, (- \cL)^{-1} (- \cA^2))\psi_S \rangle_{\pi}
\\ & = & \|\cA\psi_S \|_{\pi}^2 +  \|\cB \psi \|_{\pi}^2 \geq 0,
\end{eqnarray*}
with $\psi := (-\cL)^{-1} (- \cA^2)\psi_S$.
$\qed$
\begin{remark}
Notice also that the perturbation $ \cA^2$ is only negative semidefinite. In particular, the null space of the perturbation is (much) larger than that of the generator of the overdamped Langevin dynamics which consists of constants. The amount of improvement in the calculation of the integral in~\eqref{e:expect} using the long time average depends on the magnitude of the projection of the observable $\phi$ on the null space of $ \cA^2$. Clearly, if this projection is zero, then the inequality in~\eqref{e:var} is strict. The details of these arguments will be presented elsewhere.
\end{remark}

\section{Numerical experiments}
\begin{figure}[htb]
\centering
\begin{subfigure}[t]{0.49\textwidth} \centering
\includegraphics[width=1\linewidth]{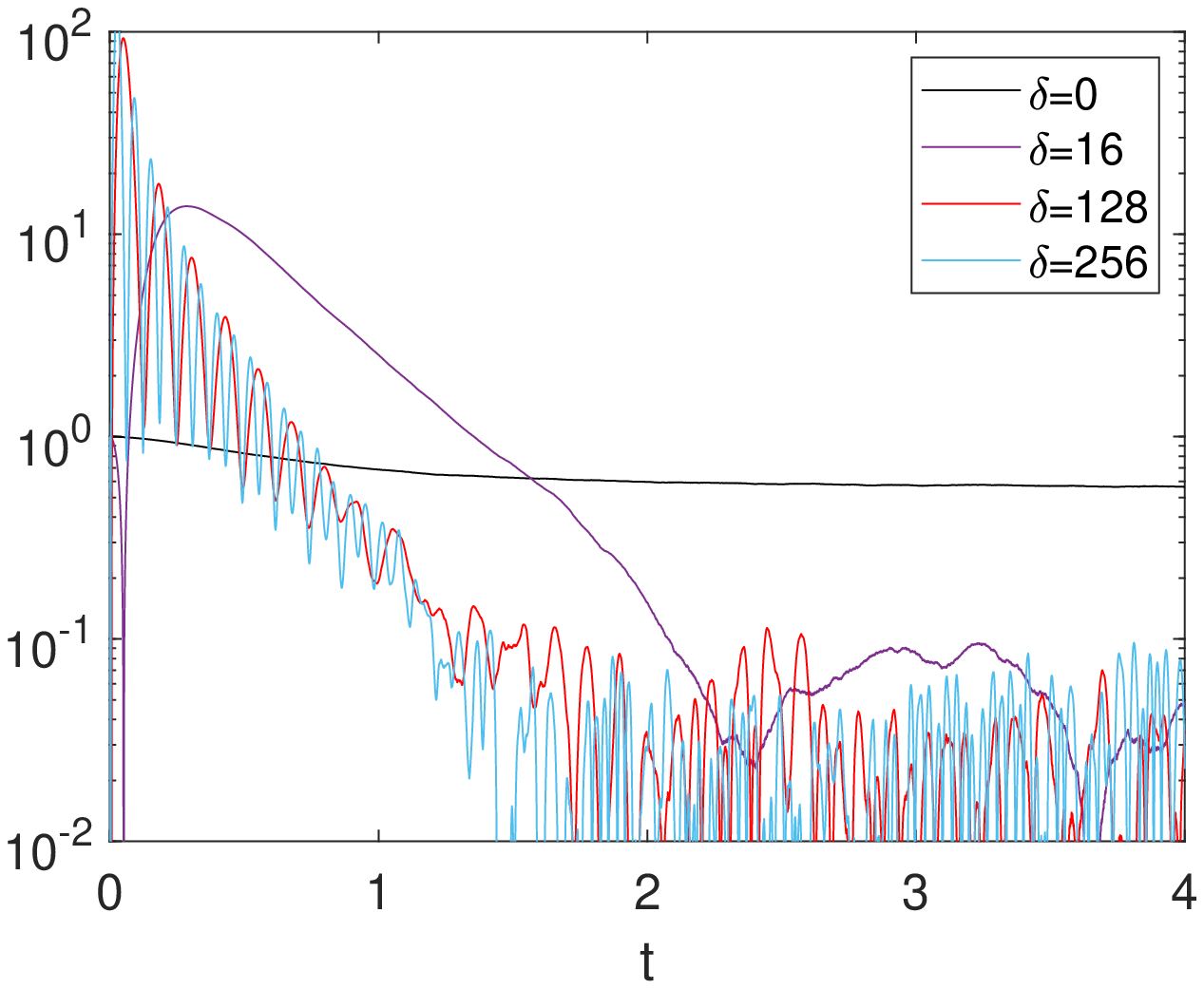}
\caption{Nonreversible Langevin dynamics~\eqref{eq:sdedet}.} \label{fig:1a}
\end{subfigure}
\begin{subfigure}[t]{0.49\textwidth} \centering
\includegraphics[width=1\linewidth]{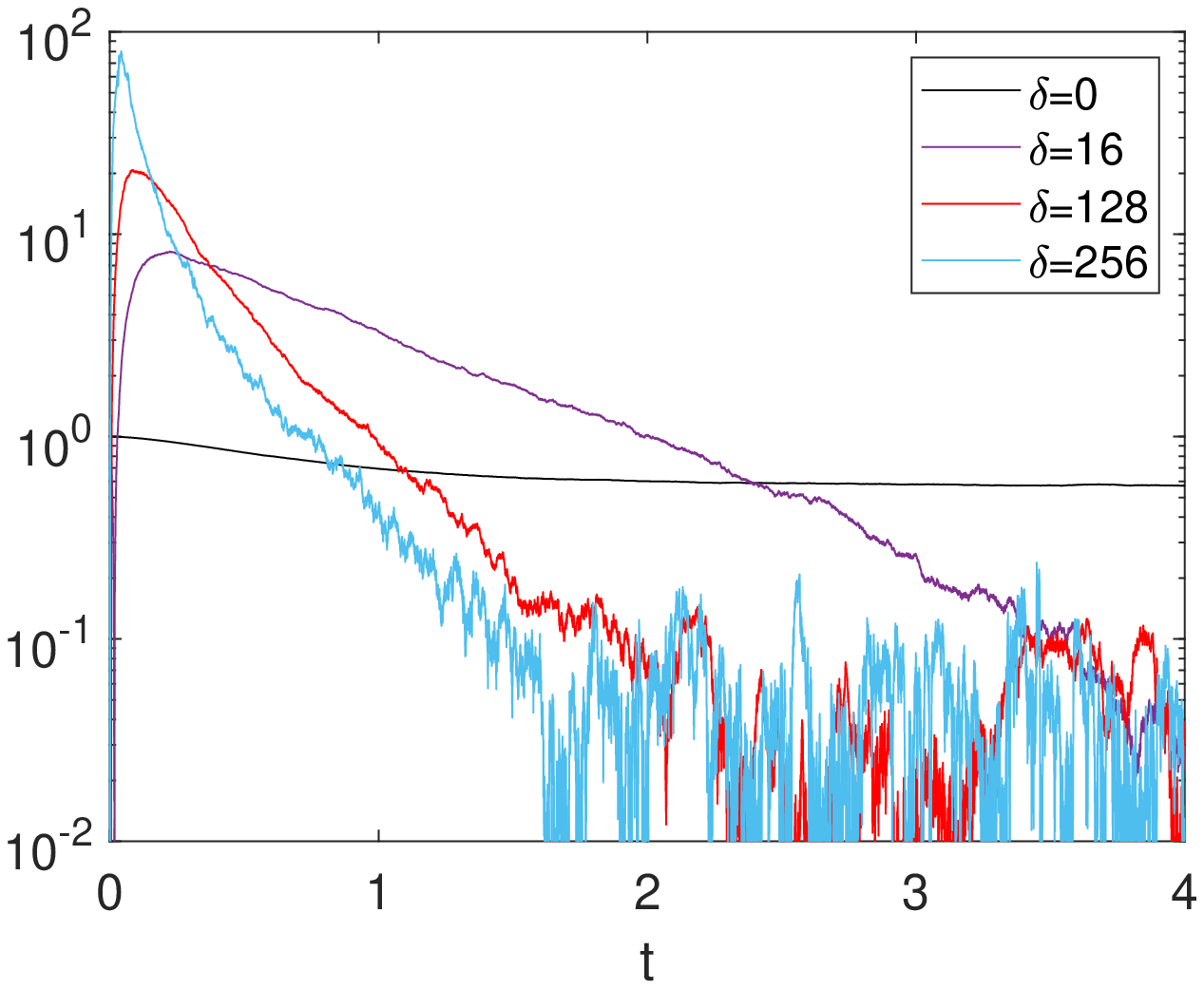}
\caption{Stratonovitch-perturbed Langevin dynamics~\eqref{eq:sdestrato}.} \label{fig:1b}
\end{subfigure}
\caption{Error evolution along time of the average over $M=10^3$ trajectories of the nonreversible and  the Stratonovitch-perturbed  Langevin dynamics  for different sizes $\delta=0,16,128,256$ of the perturbation.}
\label{fig:1}
\end{figure}

In this section, we present some numerical experiments to corroborate our theoretical findings and illustrate the features of the Stratonovitch-perturbed Langevin dynamics~\eqref{eq:sdestrato}. Although we are primarily interested in large dimensional problems, we consider for simplicity the following warped Gaussian distribution, as considered in \cite[Sec. 5.2]{DuncanLelievrePavliotis2016},
with density $\pi(x)=Z^{-1}e^{-V(x)}$ where $V(x)$ is the two-dimensional potential
potential
$V(x)=\frac{x_1^2}{100}+(x_2+bx_1^2-100b)^2,$
where the parameter $b=0.05$ is related to how warped the distribution is. For the purposes of this paper, it is sufficient to consider the family of vector fields $g(x) =  J \nabla V(x), \ J = - J^T$, for all constant skew-symmetric matrices $J$. In particular, we consider 
the vector field $g(x)$ defined by \eqref{eq:gdelta}
and we compare the effect of the nonreversible perturbation with $\theta=1$ in \eqref{eq:sdedet} (Figure \ref{fig:1a}) and the new reversible Stratonovitch perturbation  with $\theta=1/2$ in \eqref{eq:sdestrato} (Figure \ref{fig:1b}) for several sizes $\delta=1,64,256$ of the perturbation. 
We also include for reference the results for the standard overdamped Langevin equation \eqref{e:langevin}.
We consider the observable $\phi(x)=x_1^2+x_2^2$ and consider the estimator
$\frac1M\sum_{i=1}^M \phi(X^{(i)}(t)) \simeq \E(\phi(X(t))$.
We take the initial condition $X_0=(0,0)$ and we plot for $M=10^3$ independent realisations $X^{(i)}(t),i=1,\ldots,M$
the error $|\frac1M\sum_{i=1}^M \phi(X^{(i)}(t)) -\pi(\phi)|$ as a function of time $t\in[0,4]$.
The solution is approximated using the simplest Euler-Maruyama method with very small stepsize $\Delta t=10^{-5}$ (considering the It\^o formulation \eqref{eq:stratotointo}). 
We observe that although the speed of the convergence $\E(\phi(X(t))\rightarrow \pi(\phi)$ as $t\rightarrow \infty$ is very slow for the standard overdamped Langevin dynamics (see the nearly horizontal black curve for $\delta=0$), 
both perturbations lead to an increase in the speed of the convergence to equilibrium (see the transient phase for small time $t$) while reducing the asymptotic variance (see the equilibrium phase for large time $t\geq2$ where the oscillations are only due to Monte-Carlo errors), which corroborates Theorem \ref{th:pres_inv_measure} and Theorem \ref{th:reveribility}. 
In addition, the Stratonovitch perturbation yields no oscillatory behavior in contrast to the nonreversible one (see Theorem \ref{th:reversibility}). 
This feature renders the new sampling scheme more amenable to efficient numerical methods. This will be explored further in a future study.

\bibliographystyle{abbrv}
\bibliography{mybib}

\def\cprime{$'$} \def\cprime{$'$} \def\cprime{$'$} \def\cprime{$'$}
  \def\cprime{$'$} \def\cprime{$'$} \def\cprime{$'$} \def\cprime{$'$}
  \def\Rom#1{\uppercase\expandafter{\romannumeral #1}}\def\u#1{{\accent"15
  #1}}\def\Rom#1{\uppercase\expandafter{\romannumeral #1}}\def\u#1{{\accent"15
  #1}}\def\cprime{$'$} \def\cprime{$'$} \def\cprime{$'$} \def\cprime{$'$}
  \def\cprime{$'$} \def\cprime{$'$} \def\cprime{$'$}
  \def\polhk#1{\setbox0=\hbox{#1}{\ooalign{\hidewidth
  \lower1.5ex\hbox{`}\hidewidth\crcr\unhbox0}}} \def\cprime{$'$}
  \def\cprime{$'$} \def\cprime{$'$}
\begin{thebibliography}{10}

\bibitem{Andrieu2003}
C.~Andrieu, N.~de~Freitas, A.~Doucet, and M.~I. Jordan.
\newblock An introduction to {M}{C}{M}{C} for machine learning.
\newblock {\em Machine Learning}, 50(1):5--43, Jan 2003.

\bibitem{DuncanLelievrePavliotis2016}
A.~B. Duncan, T.~Leli{\`e}vre, and G.~A. Pavliotis.
\newblock Variance {R}eduction {U}sing {N}onreversible {L}angevin {S}amplers.
\newblock {\em J. Stat. Phys.}, 163(3):457--491, 2016.

\bibitem{DuncanNuskenPavliotis2017}
A.~B. Duncan, N.~N\"usken, and G.~A. Pavliotis.
\newblock Using {P}erturbed {U}nderdamped {L}angevin {D}ynamics to
  {E}fficiently {S}ample from {P}robability {D}istributions.
\newblock {\em J. Stat. Phys.}, 169(6):1098--1131, 2017.

\bibitem{DuncanPavliotisZygalakis2017}
A.~B. Duncan, G.~A. Pavliotis, and K.~C. Zygalakis.
\newblock Nonreversible {L}angevin samplers: Splitting schemes, analysis and
  implementation.
\newblock {\em ArXiv e-prints}, 2017.

\bibitem{Hwang_al1993}
C.-R. Hwang, S.-Y. Hwang-Ma, and S.-J. Sheu.
\newblock Accelerating {G}aussian diffusions.
\newblock {\em Ann. Appl. Probab.}, 3(3):897--913, 1993.

\bibitem{Hwang_al2005}
C.-R. Hwang, S.-Y. Hwang-Ma, and S.-J. Sheu.
\newblock Accelerating diffusions.
\newblock {\em Ann. Appl. Probab.}, 15(2):1433--1444, 2005.

\bibitem{LelievreNierPavliotis2013}
T.~Lelievre, F.~Nier, and G.~A. Pavliotis.
\newblock Optimal non-reversible linear drift for the convergence to
  equilibrium of a diffusion.
\newblock {\em J. Stat. Phys.}, {152}(2):{ 237--274 }, {2013}.

\bibitem{Lelievre_al_2009}
T.~Leli{\`e}vre, M.~Rousset, and G.~Stoltz.
\newblock {\em Free energy computations}.
\newblock Imperial College Press, London, 2010.
\newblock A mathematical perspective.

\bibitem{Mattingly10con}
J.~C. Mattingly, A.~M. Stuart, and M.~V. Tretyakov.
\newblock Convergence of numerical time-averaging and stationary measures via
  {P}oisson equations.
\newblock {\em SIAM J. Numer. Anal.}, 48(2):552--577, 2010.

\bibitem{Pavl2014}
G.~A. Pavliotis.
\newblock {\em Stochastic processes and applications}, volume~60 of {\em Texts
  in Applied Mathematics}.
\newblock Springer, New York, 2014.
\newblock Diffusion processes, the Fokker-Planck and Langevin equations.

\bibitem{Bellet_Spiliopoulos_2015}
L.~Rey-Bellet and K.~Spiliopoulos.
\newblock Irreversible langevin samplers and variance reduction: a large
  deviations approach.
\newblock {\em Nonlinearity}, 28(7):2081, 2015.

\bibitem{Bellet_Spiliopoulos_2016}
L.~Rey-Bellet and K.~Spiliopoulos.
\newblock Improving the convergence of reversible samplers.
\newblock {\em J. Stat. Phys.}, 164(3):472--494, 2016.

\bibitem{Stuart2010}
A.~M. Stuart.
\newblock Inverse problems: a {B}ayesian perspective.
\newblock {\em Acta Numer.}, 19:451--559, 2010.

\end{thebibliography}

\end{document}